\documentclass[11pt]{article}
\usepackage{amsmath,amsthm,amsfonts}
\newcommand\pd[2]{\frac{\partial#1}{\partial#2}}
\renewcommand{\=}{\doteq}
\newcommand{\p}{\partial}
\newtheorem{thm}{Theorem}[section]

 \newtheorem{lemma}[thm]{Lemma}
\theoremstyle{definition}

\theoremstyle{definition}
 
\numberwithin{equation}{section}
\numberwithin{equation}{section}
\begin{document}
\title{\bf  Moufang symmetry III.\\
Integrability of generalized Lie equations
}
\author{Eugen Paal}
\date{}
\maketitle
\thispagestyle{empty}
\begin{abstract}
Integrability of generalized Lie equations of a local analytic Moufang loop is inquired. 
\end{abstract}

\section{Introduction}

In this paper we proceed explaing the Moufang symmetry. The paper can be seen as a continuation of \cite{Paal1,Paal2}.

\section{Generalized Lie equations}

In \cite{Paal1} the \emph{generalized Lie equations} (GLE) of a local analytic Moufang loop $G$ were found. These read
\begin{subequations}
\label{gle}
\begin{align}
w^{s}_{j}(g)\pd{(gh)^{i}}{g^{s}}+u^{s}_{j}(h)\pd{(gh)^{i}}{h^{s}}+u^{i}_{j}(gh)&=0\\
v^{s}_{j}(g)\pd{(gh)^{i}}{g^{s}}+w^{s}_{j}(h)\pd{(gh)^{i}}{h^{s}}+v^{i}_{j}(gh)&=0\\
u^{s}_{j}(g)\pd{(gh)^{i}}{g^{s}}+v^{s}_{j}(h)\pd{(gh)^{i}}{h^{s}}+w^{i}_{j}(gh)&=0
\end{align}
\end{subequations}
where $gh$ is the product of $g$ and $h$, and the auxiliary functions $u^s_j$, $v^s_j$ and $w^s_j$ are related with the  constraint
\begin{equation}
u^s_j(g)+v^s_j(g)+w^s_j(g)=0
\end{equation}
In this paper we inquire integrability of GLE (\ref{gle}a--c). Triality \cite{Paal2} considerations are wery helpful.

\section{Generalized Maurer-Cartan equations and Yamagutian}

Recall from \cite{Paal1} that for $x$ in  $T_e(G)$ the infinitesimal translations of $G$ are defined by
\begin{equation*}
L_x\=x^j u^s_j(g)\pd{}{g^s},\quad
R_x\=x^j v^s_j(g)\pd{}{g^s},\quad
M_x\=x^j w^s_j(g)\pd{}{g^s}\quad \in T_g(G)
\end{equation*}
with constriant
\begin{equation*}
L_x+R_x+M_x=0
\end{equation*}
Following triality \cite{Paal2} define the Yamagutian $Y(x;y)$ by
\begin{equation*}
6Y(x;y)=[L_x,L_y]+[R_x,R_y]+[M_x,M_y]
\end{equation*}
We know  from \cite{Paal2} the generalized Maurer-Cartan equations:
\begin{subequations}
\label{m-c}
\begin{align} 
[L_{x},L_{y}]&=L_{[x,y]}-2[L_{x},R_{y}]\\
[R_{x},R_{y}]&=R_{[y,x]}-2[R_{x},L_{y}]\\
[L_{x},R_{y}]&=[R_{x},L_{y}],\quad \forall x,y \in T_e(G) 
\end{align}
\end{subequations}
The latter can be written \cite{Paal2} as follows:
\begin{subequations}
\label{lr-yam}
\begin{align}
[L_{x},L_{y}]=2Y(x;y)+\frac{1}{3}L_{[x,y]}+\frac{2}{3}R_{[x,y]}\\
[L_{x},R_{y}]=-Y(x;y)+\frac{1}{3}L_{[x,y]}-\frac{1}{3}R_{[x,y]}\\
[R_{x},R_{y}]=2Y(x;y)-\frac{2}{3}L_{[x,y]}-\frac{1}{3}R_{[x,y]}
\end{align}
\end{subequations}
Define the  (secondary) auxiliary functions of $G$ by
\begin{align*}
u^s_{jk}(g)
&\=u^p_k(g)\pd{u^s_j(g)}{g^p}-u^p_j(g)\pd{u^s_k(g)}{g^p}\\
v^s_{jk}(g)
&\=v^p_k(g)\pd{v^s_j(g)}{g^p}-v^p_j(g)\pd{v^s_k(g)}{g^p}\\
w^s_{jk}(g)
&\=w^p_k(g)\pd{w^s_j(g)}{g^p}-w^p_j(g)\pd{w^s_k(g)}{g^p}
\end{align*}
The Yamaguti functions $Y^i_{jk}$ are defined by
\begin{equation*}
6Y^s_{jk}(g)\=u^s_{jk}(g)+v^s_{jk}(g)+w^s_{jk}(g)
\end{equation*}
Evidently,
\begin{align*}
[L_x,L_y]
&=-x^jy^k u^s_{jk}(g)\pd{}{g^s}\\
[R_x,R_y]
&=-x^jy^k v^s_{jk}(g)\pd{}{g^s}\\
[M_x,M_y]
&=-x^jy^k w^s_{jk}(g)\pd{}{g^s}
\end{align*}
By adding the above formulae, we get
\begin{equation*}
Y(x;y)=-x^jy^k Y^s_{jk}(g)\pd{}{g^s}
\end{equation*}

\begin{lemma}
One has
\begin{subequations}
\label{uvw}
\begin{align}
u^i_{jk}
&\=2Y^i_{jk}+\frac{1}{3}C^s_{jk}(u^i_ṣ+2v^i_ṣ)\\
v^i_{jk}
&\=2Y^i_{jk}-\frac{1}{3}C^s_{jk}(2u^i_ṣ+v^i_ṣ)\\
w^i_{jk}
&\=2Y^i_{jk}+\frac{1}{3}C^s_{jk}(u^i_ṣ-v^i_ṣ)
\end{align}
\end{subequations}
\end{lemma}

\begin{proof}
To see (\ref{uvw}a,b) use (\ref{lr-yam}a,c) . To see (\ref{uvw}c) calculate by using (\ref{lr-yam}):
\begin{align*}
[M_x,M_y]
&=[L_x+R_x,L_y+R_y]\\
&=[L_x,L_y]+[L_x,R_y]+[R_x,L_y]+[R_x,R_y]\\
&=[L_x,L_y]+2[L_x,R_y]+[R_x,R_y]\\
&=2Y(x;y)+\frac{1}{3}\left(L_{[x,y]}-R_{[x,y]}\right)
\end{align*}
and  (\ref{uvw}b) easily follows.
\end{proof}

\section{Integrability conditions}

\begin{thm}
The integrability conditons of the GLE (\ref{gle}a--c) read
\begin{equation}
\label{gle2yam}
Y^s_{jk}(g)\pd{(gh)^i}{g^s}+Y^s_{jk}(h)\pd{(gh)^i}{g^s}=Y^i_{jk}(gh)
\end{equation}
\end{thm}

\begin{proof}
We differentiate the GLE and use
\begin{equation}
\label{int}
\frac{\p^{2}(gh)^i}{\p g^{j}\p g^{k}}=\frac{\p^{2}(gh)^i}{\p g^{k}\p g^{j}},\quad
\frac{\p^{2}(gh)^i}{\p g^{j}\p h^{k}}=\frac{\p^{2}(gh)^i}{\p g^{k}\p h^{j}},\quad
\frac{\p^{2}(gh)^i}{\p h^{j}\p h^{k}}=\frac{\p^{2}(gh)^i}{\p h^{k}\p h^{j}}
\end{equation}
First differentiate (\ref{gle}a) with respect to $g^p$ and $h^p$:
\begin{subequations}
\label{int1}
\begin{align}
\pd{w^s_j(g)}{g^p}\pd{(gh)^i}{g^s}
+w^s_j(g)\frac{\p^2(gh)^i}{\p g^p \p g^s}
+u^s_j(h)\frac{\p^2(gh)^i}{\p g^p \p h^s}
&=
-\pd{u^i_j(gh)}{(gh)^s}\pd{(gh)^s}{g^p}\\
w^s_j(g)\frac{\p^2(gh)^i}{\p h^p \p g^s}
+\pd{u^s_j(g)}{h^p}\pd{(gh)^i}{h^s}
+u^s_j(h)\frac{\p^2(gh)^i}{\p h^p \p h^s}
&=
-\pd{u^i_j(gh)}{(gh)^s}\pd{(gh)^s}{h^p}
\end{align}
\end{subequations}
Now multiply (\ref{int1}a) by $w^p_k(g)$ and (\ref{int1}b) by $u^p_k(g)$ and add the resulting formulae. On the right hand side of the resulting formula use again the GLE (\ref{gle}a); then transpose the indexes $j$ and $k$ and subtract the result from the previous one. Then it turns out that due to (\ref{int}) all terms with the second order partial derivatives vanish and result reads
\begin{equation}
\label{gle2a}
w^{s}_{jk}(g)\pd{(gh)^{i}}{g^{s}}+u^{s}_{jk}(h)\pd{(gh)^{i}}{h^{s}}=u^{i}_{jk}(gh)
\end{equation}
By acting analogously with GLE (\ref{gle}b,c) we get
\begin{subequations}
\label{gle2b}
\begin{align}
v^{s}_{jk}(g)\pd{(gh)^{i}}{g^{s}}+w^{s}_{jk}(h)\pd{(gh)^{i}}{h^{s}}&=v^{i}_{jk}(gh)\\
u^{s}_{jk}(g)\pd{(gh)^{i}}{g^{s}}+v^{s}_{jk}(h)\pd{(gh)^{i}}{h^{s}}&=w^{i}_{jk}(gh)
\end{align}
\end{subequations}
Now add (\ref{gle2a}), (\ref{gle2b}a) and (\ref{gle2b}b) to obtain (\ref{gle2yam}).
 
It remains to show that  (\ref{gle2a}), (\ref{gle2b}a) and (\ref{gle2b}b)  are equivalent to (\ref{gle2yam}).
By using (\ref{uvw}a--c) calculate
\begin{align*}
w^{s}_{jk}(g)\pd{(gh)^{i}}{g^{s}}+u^{s}_{jk}(h)\pd{(gh)^{i}}{h^{s}}-u^{i}_{jk}(gh)
&\overset{(\ref{uvw}\text{a,c})}{=}\\
v^{s}_{jk}(g)\pd{(gh)^{i}}{g^{s}}+w^{s}_{jk}(h)\pd{(gh)^{i}}{h^{s}}-v^{i}_{jk}(gh)
&\overset{(\ref{uvw}\text{b,c})}{=}\\
u^{s}_{jk}(g)\pd{(gh)^{i}}{g^{s}}+v^{s}_{jk}(h)\pd{(gh)^{i}}{h^{s}}-w^{i}_{jk}(gh)
&\overset{(\ref{uvw}\text{a--c})}{=}\\
&=2\left(Y^s_{jk}(g)\pd{(gh)^i}{g^s}+Y^s_{jk}(h)\pd{(gh)^i}{g^s}-Y^i_{jk}(gh)\right)
\tag*{\qed}
\end{align*}
\renewcommand{\qed}{}
\end{proof}

\section*{Acknowledgement}

Research was in part supported by the Estonian Science Foundation, Grant 6912.

\bigskip\noindent
Department of Mathematics\\
Tallinn University of Technology\\
Ehitajate tee 5, 19086 Tallinn, Estonia\\ 
E-mail: eugen.paal@ttu.ee

\end{document}